\documentclass[twoside]{amsart}
\usepackage[centertags]{amsmath}

\usepackage{amsfonts}
\usepackage{amssymb}
\usepackage{amsthm}
\usepackage{eucal}
\usepackage[cmtip,all]{xy}
\usepackage{graphicx}
\newtheorem{thm}[equation]{Theorem}

\newtheorem{prop}[equation]{Proposition}
\theoremstyle{definition}
\newtheorem{defn}[equation]{Definition}
\theoremstyle{remark}
\newtheorem{rem}[equation]{Remark}
\newtheorem{exm}[equation]{Example}
\numberwithin{equation}{section}

\newcommand{\set}[1]{\left\{#1\right\}}

\newcommand{\To}{\longrightarrow}

\newcommand{\C}[1]{\mathcal{#1}}

\newcommand{\EZDIAG}[5]{\xymatrix@C+=2.5cm{*+[r]{#1}
\ar@(u,l)_(0.62){\displaystyle #5}[]
\ar@<.5ex>^-{#3}[r]&\ar@<.5ex>^-{#4}[l]#2}}

\def\op{{op}}

\def\r{\rightarrow} 
\def\l{\leftarrow} 

\def\hom{\operatorname{Hom}}

\def\ho{\operatorname{Ho}}

\def\ext{\operatorname{Ext}}
\def\trace{\operatorname{trace}}

\def\st{\stackrel} 

\def\coker{\operatorname{Coker}}
\renewcommand{\ker}{\operatorname{Ker}}

\def\Z{\mathbb{Z}}

\def\S{\Sigma}

\newcommand{\grupo}[1]{\langle #1\rangle}

\numberwithin{equation}{section}

\begin{document}

\title{A triangulated category without models}%
\author{Fernando Muro}%
\address{Universitat de Barcelona, Departament d'$\grave{\text{A}}$lgebra i Geometria, Gran Via de les Corts Catalanes, 585, 08007 Barcelona, Spain}
\email{fmuro@ub.edu}

\subjclass{18E30, 55P42, 16E40}
\keywords{Triangulated category, stable model category, cohomology of categories, Toda bracket, Mac Lane cohomology}%

\begin{abstract}
We exhibit a triangulated category which is neither the stable category of a Frobenius category nor a full triangulated subcategory of the homotopy category of a stable model category.
\end{abstract} \maketitle

\section{Introduction}

Triangulated categories are fundamental tools in both algebra and topology. In algebra they use to arise as the stable category of a Frobenius category (\cite[4.4]{shc}, \cite[IV.3 Exercise 8]{mha}). In topology they usually appear as a full triangulated subcategory of the homotopy category of a Quillen stable model category (\cite[7.1]{hmc}). The triangulated categories which belong, up to exact equivalence, to one of these two families will be termed \emph{algebraic} and \emph{topological}, respectively. We borrow this terminology from \cite[3.6]{odgc} and \cite{avttc}. Algebraic triangulated categories are generally also topological, but there are many well-known examples of topological triangulated categories which are not algebraic. 

In the present paper we produce the first example of a triangulated category which is neither algebraic nor topological.

Let $\C{P}(R)$ be the category of finitely generated projective  modules over a commutative ring $R$.

\begin{thm}\label{rara}
The pair given by the category $\C{P}(\Z/4)$ together with the identity functor can be endowed with a unique triangulated category structure. 
\end{thm}

This theorem is proved in Section \ref{ttc}. There we also show that
\begin{equation}\label {rari}
\Z/4\st{2}\To\Z/4\st{2}\To\Z/4\st{2}\To\Z/4,
\end{equation}
is an exact triangle in $\C{P}(\Z/4)$, therefore $\C{P}(\Z/4)$ cannot be algebraic. Indeed given an object $X$ in an algebraic triangulated category $\C{T}$ and an exact triangle
$$X\st{2\cdot 1_X}\To X\To Y\To\S X,$$
the equation $2\cdot 1_Y=0$ holds, compare \cite[3.6]{odgc} and \cite{avttc}. This needs not happen in a topological triangulated category. This fact is illustrated by the classical example of $X=S$ the sphere spectrum in the Spanier-Whitehead category. Nevertheless we prove as a main theorem that the triangulated category $\C{P}(\Z/4)$ is not topological.


\begin{thm}\label{raro}
The triangulated category $\C{P}(\Z/4)$ with the identity translation functor is not exact equivalent to a full triangulated subcategory of the homotopy category $\ho\C{M}$ of a stable model category $\C{M}$. 
\end{thm}

The proof of this theorem is completed at the very end of the paper. Unlike in the algebraic case, one cannot prove that $\C{P}(\Z/4)$ is not topological by using local properties. In order to prove Theorem \ref{raro} we use the existence of the \emph{universal Toda bracket} in the cohomology of the homotopy category
$$\grupo{\C{M}}\in H^3(\ho\C{M},\hom_{\ho\C{M}}(\Sigma,-)),$$
as defined in \cite{ccglg}. This universal Toda bracket determines all Toda brackets in $\ho\C{M}$, and Toda brackets determine the exact triangles, see \cite[Theorem 13.2]{shc}. We will explicitly compute $H^3(\C{P}(\Z/4),\hom_{\Z/4})$ and show that no element in this group would yield (\ref{rari}) as an exact triangle.

The cohomological theory of triangulated categories is studied in \cite{ccctc1,ccctc2}. There, given an additive category $\C{A}$, a self equivalence $\Sigma\colon\C{A}\st{\sim}\r\C{A}$, and a cohomology class
$$\nabla\in H^3(\C{A},\hom_\C{A}(\Sigma,-)),$$
we describe conditions on $\nabla$, or rather on a lift of $\nabla$ to the translation cohomology group $H^3(\C{A},\Sigma)$, under which the family of exact triangles determined by $\nabla$ yields a triangulated structure on $\C{A}$ with translation functor $\Sigma$. Not only algebraic and topological triangulated categories arise in this way. More generally, the triangulated homotopy category of a stable $S$-category $\C{S}$, in the sense of \cite[7]{ktsc}, is associated to the classifying cohomology class (\cite{catc}) of the groupoid-enriched category obtained by taking fundamental groupoid on the morphism simplicial sets of $\C{S}$. This cohomology class can be regarded as the first Postnikov invariant of $\C{S}$, compare \cite{otsc}. The theory of stable $\infty$-categories developed in \cite{dag1} is equivalent to the theory of stable $S$-categories, compare \cite{qcvsc}, so the homotopy category of a stable $\infty$-category is also cohomologically triangulated in the sense of \cite{ccctc1,ccctc2}. The results of this paper show that $\C{P}(\Z/4)$ is not cohomologically triangulated so it cannot be obtained by any of the procedures described above.


\subsection*{Acknowledgements}

Stefan Schwede's talk on algebraic vs. topological triangulated categories at the ICM 2006 Satellite Workshop on Triangulated Categories, held in Leeds, stimulated my interest on this subject. I am grateful to Bernhard Keller for conversations on the results of this paper. I feel indebted 
to Amnon Neeman, who suggested the possibility of constructing a triangulated structure on $\C{P}(\Z/4)$ by using Heller's theory.



\section{The triangulated category}\label{ttc}

Recall that a \emph{triangulated category} is an additive category $\C{T}$ together with a self-equivalence $\S\colon
\C{T}\st{\sim}\r\C{T}$, called \emph{translation functor}, and a collection of diagrams
\begin{equation}\label{tri}
A\st{f}\To B\st{i}\To C\st{q}\To \S A,
\end{equation}
called \emph{distinguished} or \emph{exact triangles}, which satisfy Puppe's axioms (\cite{fssht}) and Verdier's octahedral axiom
(\cite{cd}), see \cite{triang} for a unified reference. If we do not require the octahedral axiom then we speak of
a \emph{pretriangulated category}, as considered for example in \cite{shc}. An \emph{exact functor} between (pre)triangulated categories is an additive functor commuting with the translation functors up to natural isomorphism and preserving exact triangles.

Heller developed in \cite{shc} a method to ``count'' pretriangulated structures on an additive category $\C{T}$
which is the category of injective objects in a Frobenius category $\C{B}$, with a given translation functor
$\S$. The following proposition is essentially a restatement of \cite[Theorem 16.4]{shc} for $\C{B}$ the category of finitely
generated $\Z/4$-modules.

\begin{prop}\label{hel}
There is a bijective correspondence between the pretriangulated category structures on $\C{P}(\Z/4)$ with
identity translation functor and the automorphisms of the identity functor in $\C{P}(\Z/2)$.
\end{prop}

\begin{proof}
The category $\C{B}$ of finitely generated $\Z/4$-modules is a Frobenius category, see for instance
\cite[Theorems 4.35
and 4.37]{ihh}. Any object in $\C{B}$ is isomorphic to a finite direct sum of copies of $\Z/2$ and $\Z/4$. Moreover, $\Z/4$
is an injective hull of $\Z/2$, therefore the stable category $\C{B}^\square$ is equivalent to $\C{P}(\Z/2)$ through the functor
$\C{P}(\Z/2)\subset\C{B}\twoheadrightarrow \C{B}^\square$. The suspension functor $S$ in $\C{B}^\square$ is
determined by the choice of short exact sequences in $\C{B}$
\begin{equation*}
X\st{j}\hookrightarrow CX\st{r}\twoheadrightarrow SX,
\end{equation*}
with $CX$ injective, for each $X$ in $\C{B}$. For our convenience we choose $CX$ to be an injective hull of $X$, hence $SX$ is always a $\Z/2$-vector space. For $X$ injective we can take $CX=X$ and $SX=0$, and for $X$ a $\Z/2$-vector space $SX=X$ and $jr=2\colon CX\r CX$. Notice that $S$ in $\C{B}^\square$ restricts to the identity functor in $\C{P}(\Z/2)$. Notice also that any automorphism $\phi$ of the identity functor in $\C{P}(\Z/2)$ satisfies $\phi+\phi=0$. Now the
proposition follows from \cite[Theorem 16.4]{shc}.
\end{proof}

The category $\C{P}(\Z/4)$ has at least
the pretriangulated structure associated to the identity automorphism of the identity functor in $\C{P}(\Z/2)$, 
for which we will prove in this section the octahedral axiom. By Proposition \ref{solo1} this is the only
automorphism of the identity functor in $\C{P}(\Z/2)$. This will complete the proof of Theorem \ref{rara}. 

We now recall some useful definitions from \cite{triang}. A \emph{candidate triangle} in a pretriangulated category is a diagram as (\ref{tri}) where $if$, $qi$, and
$(\S f)q$ are zero morphisms. 
A \emph{morphism} of candidate triangles is a commutative diagram
\begin{equation}\label{candmor}
\xymatrix{A\ar[r]^f\ar[d]^{k_0}&B\ar[r]^i\ar[d]^{k_1}&C\ar[r]^q\ar[d]^{k_2}&\Sigma A\ar[d]^{\Sigma k_0}\\
A'\ar[r]^{f'}&B'\ar[r]^{i'}&C'\ar[r]^{q'}&\Sigma A'} 
\end{equation}
This morphism $k=(k_0,k_1,k_2)$ is \emph{homotopic} to another morphism $k'=(k_0',k_1',k_2')$ between the same candidate triangles if there are morphisms
\begin{eqnarray*}
B&\st{\alpha_1}\To& A',\\
C&\st{\alpha_2}\To &B',\\
\S A&\st{\alpha_0}\To& C',
\end{eqnarray*}
such that
\begin{eqnarray*}
-k_1+k_1'&=&f'\alpha_1+\alpha_2i,\\
-k_2+k_2'&=&i'\alpha_2+\alpha_0q,\\
-\Sigma k_0+\Sigma k_0'&=&q'\alpha_0+\Sigma (\alpha_1f).
\end{eqnarray*}
A candidate triangle is \emph{contractible} if the identity morphism is homotopic to the zero morphism.

We also recall that the \emph{mapping cone} of a morphism like (\ref{candmor}) is the candidate triangle
\begin{equation*}
B\oplus A'\st{\begin{scriptsize}
\left(\begin{array}{cc}
-i&0\\
k_1&f'
\end{array}\right)\end{scriptsize}
}\To C\oplus B'\st{\begin{scriptsize}
\left(\begin{array}{cc}
-q&0\\
k_2&i'
\end{array}\right)
              \end{scriptsize}
}\To \Sigma A\oplus C'\st{\begin{scriptsize}
\left(\begin{array}{cc}
-\Sigma f&0\\
\Sigma k_0&q'
\end{array}\right)\end{scriptsize}}\To \Sigma B\oplus\Sigma A'.
\end{equation*}
The mapping cones of two homotopic morphisms are isomorphic.

\begin{prop}\label{lose}
Consider the pretriangulated category structure on $\C{P}(\Z/4)$ with identity translation functor associated to the identity automorphism of the identity functor in $\C{P}(\Z/2)$. The exact triangles of this triangulated category are the candidate triangles which are isomorphic to the direct sum 
of a contractible candidate triangle with a candidate triangle of the form
$$X\st{2}\To X\st{2}\To X\st{2}\To X.$$
\end{prop}

\begin{proof}
Contractible candidate triangles are always exact, see \cite[Proposition 1.3.8]{triang}. Any object $X$ in
$\C{P}(\Z/4)$ is, up to isomorphism, of the form $CY$ for some $\Z/2$-vector space $Y$. The canonical resolution of
$Y$ given by the choice of cones and suspensions in the proof of Proposition \ref{hel} takes the form
\begin{equation*}\tag{a}
CY\st{2}\To CY\st{2}\To CY\st{2}\To CY\r\cdots,
\end{equation*}
therefore the description of exact triangles in $\C{P}(\Z/4)$ given in \cite{shc} shows that the triangle of the
statement is exact.

Let us now prove that all exact triangles can be decomposed as indicated in the statement. Let
\begin{equation*}\tag{b}
A\st{f}\To B\st{i}\To C\st{q}\To A
\end{equation*}
be an exact triangle. The kernel of $f$ in the category of $\Z/4$-modules decomposes as a direct sum $\ker f\cong X\oplus Y$ with $X$ injective and $Y$ a $\Z/2$-vector space. Since (a) is a minimal injective resolution of $Y$, combining elementary homological algebra with axiom (Tr3) for pretriangulated categories, see \cite{triang}, we obtain a morphism $k$ of candidate triangles given by monomorphisms
\begin{equation*}
\xymatrix{Y\ar@{^{(}-->}[d]\ar@{^{(}-->}[r]^j&CY\ar[r]^2\ar@{^{(}->}[d]^{k_0}&CY\ar[r]^2\ar@{^{(}->}[d]^{k_1}&CY\ar[r]^2\ar@{^{(}->}[d]^{k_2}&CY\ar@{^{(}->}[d]^{k_0}\\
\ker f\ar@{^{(}-->}[r]&A\ar[r]^{f}&B\ar[r]^{i}&C\ar[r]^{q}&A} 
\end{equation*}
such that the square with dashed arrows on the left (which is not part of the morphism $k$) commutes.
Since $CY$ is injective the exact triangle (b) decomposes as the direct sum of the following two candidate triangles 
\begin{equation}\tag{c}
CY\st{2}\To CY\st{2}\To CY\st{2}\To CY,
\end{equation}
\begin{equation}\tag{d}
\coker k_0\st{\bar{f}}\To \coker k_1\st{\bar{i}}\To \coker k_2\st{\bar{q}}\To \coker k_0.
\end{equation}
The axioms of a pretriangulated category show that if we extend (b) to a $3$-periodic cochain complex then the resulting cochain complex is acyclic, as it happens with (c), therefore the long exact sequence in cohomology yields the same property for (d). In particular we can decompose (d) as a diagram
\begin{equation*}
X_0\st{i_0}\hookrightarrow\coker k_0\st{p_0}\twoheadrightarrow X_1\st{i_1}\hookrightarrow\coker k_1\st{p_1}\twoheadrightarrow X_2 \st{i_2}\hookrightarrow \coker k_2\st{p_2}\twoheadrightarrow X_0\st{i_0}\hookrightarrow \coker k_0.
\end{equation*}
where each subdiagram $\bullet \hookrightarrow\bullet \twoheadrightarrow\bullet$ is a short exact sequence.
By construction $X_0\cong X$, the injective direct summand of $\ker f$. Since proyectives and injectives coincide we derive that all $X_i$ are injective, as the $\coker k_i$, so we can take morphisms
\begin{equation*}
X_0\st{r_0}\twoheadleftarrow\coker k_0\st{l_0}\hookleftarrow X_1\st{r_1}\twoheadleftarrow\coker k_1\st{l_1}\hookleftarrow X_2 \st{r_2}\twoheadleftarrow \coker k_2\st{l_2}\hookleftarrow X_0\st{r_0}\twoheadleftarrow \coker k_0
\end{equation*}
such that $r_si_s=1_{X_s}$, $p_sl_s=1_{X_{s+1}}$, and $i_sr_s+l_sp_s=1_{\coker k_s}$ for all $s\in\Z/3=\set{0,1,2}$. Now one can easily check that the morphisms $\alpha_s=r_sl_{s-1}$ yield a contracting homotopy for the candidate triangle (d).
\end{proof}

Now we use this explict description of the exact triangles in $\C{P}(\Z/4)$ to prove the octahedral axiom.

\begin{prop}
The pretriangulated category structure on $\C{P}(\Z/4)$ with identity translation functor associated to the identity automorphism of the identity functor in $\C{P}(\Z/2)$ is indeed triangulated.
\end{prop}

\begin{proof}
We are going to show that $\C{P}(\Z/4)$ satisfies axiom (Tr4') in \cite[Definition 1.3.13]{triang}, which is equivalent to Verdier's octahedral axiom. We have to show that for any morphism of exact triangles $k\colon T\r T'$ we can modify $k_2$ so that the mapping cone of $k$ is again an exact triangle. The mapping cone of a morphism only depends on its homotopy class, and any morphism with contractible source or target is homotopic to the zero morphism, therefore by Proposition \ref{lose} we only need to consider the case where $T$ and $T'$ are defined by two objects $X$ and $X'$ in $\C{P}(\Z/4)$ and the multiplication by $2$, as in the statement of Proposition \ref{lose}. The exact triangle given by an object $Y$ and multiplication by $2$ will be denoted by $Y(2)$, so $T=X(2)$ and $T'=X'(2)$.

We can suppose without loss of generality that 
\begin{eqnarray*}
X&=&X_1\oplus X_2\oplus X_3,\\ 
X'&=&X_1\oplus X_2\oplus X_4,\\
k_0&=&\left(\begin{array}{ccc}
1&0&0\\
0&2&0\\
0&0&0
\end{array}\right).
\end{eqnarray*} 
Let 
\begin{eqnarray*}
k_2'&=&k_1+\left(\begin{array}{ccc}
0&0&0\\
0&2&0\\
0&0&0
\end{array}\right).
\end{eqnarray*}
Obviously $k'=(k_0,k_1,k_2')\colon T\r T'$ is also a morphism of candidate triangles. We are going to show that the mapping cone of $k'$ is exact. Since $k$ is a morphism $2k_0=2k_1$, therefore there exists $g\colon X\r X'$ such that $k_1=k_0+2g$. Then $k'$ is homotopic to $k''=(k_0,k_0,k_2'')$, where
\begin{eqnarray*}
k_2''&=&\left(\begin{array}{ccc}
1&0&0\\
0&0&0\\
0&0&0
\end{array}\right),
\end{eqnarray*}
so it will suffice to show that the mapping cone of $k''$ is exact. The mapping cone of $k''$ is isomorphic to the direct sum of four candidate triangles, namely $X_3(2)$, $X_4(2)$, the mapping cone of the identity morphism in $X_1(2)$ (which is contractible), and
the mapping cone of $(2,2,0)\colon X_2(2)\r X_2(2)$, which we denote by $T_2$. The candidate triangle $T_2$ is also exact since we have an isomorphism
$$\left(\left(\begin{array}{cc}
1&0\\
1&1
\end{array}
\right),1,1\right)\colon T_2\st{\cong}\To (X_2\oplus X_2)(2).$$
Now this proposition follows from Proposition \ref{lose}.
\end{proof}

\begin{rem}
One can actually check that given a prime $p$ and a possitive integer $n$ the category $\C{P}(\Z/p^n)$ with the identity translation functor has a pretriangulated structure if and only if $n=1$ or $p=n=2$. In all cases the pretriangulated structure is unique and triangulated. For $n=1$, $\C{P}(\Z/p)$ is the stable category of the Frobenius abelian category of finitely generated $\Z/p^2$-modules, so $\C{P}(\Z/p)$ is algebraic. This illustrates the singularity of the triangulation of $\C{P}(\Z/4)$.
\end{rem}

\section{Cohomology of categories and Mac Lane cohomology}

Let us recall from \cite{csc} the definition of the Baues-Wirsching cohomology of categories.

\begin{defn}\label{cc}
Let $\C{C}$ be a category. A \emph{$\C{C}$-bimodule} $L$ is a functor
$L\colon\C{C}^\op\times\C{C}\r\mathbf{Ab}$, where $\mathbf{Ab}$ is the category of abelian groups.
The \emph{Baues-Wirsching
complex} $F^*(\C{C},L)$ is given by the following products
indexed by all sequences of morphisms of length $n$ in $\C{C}$
\begin{eqnarray*}
F^n(\C{C},L)&=&\prod_{X_0\st{\sigma_1}\l\cdots\st{\sigma_n}\l X_n}L(X_n,X_0),\;\;n\geq0.
\end{eqnarray*}
In this formula we assume that a sequence of length $0$ is an
object $X_0$ in $\C{C}$ which we also identify with the identity
morphism $1_{X_0}$. The coordinate of $c\in F^n(\C{C},L)$ in
$X_0\st{\sigma_1}\l\cdots\st{\sigma_n}\l X_n$ will be
denoted by $c(\sigma_1,\dots,\sigma_n)$. The value of the differential $d$ over an $n$-cochain $c$, $n\geq 1$, is
defined as
\begin{eqnarray*}
  d(c)(\sigma_1,\dots,\sigma_{n+1}) &=& L(X_{n+1},\sigma_1)c(\sigma_2,\dots,\sigma_{n+1}) \\
   && +\sum_{i=1}^{n}(-1)^ic(\sigma_1,\dots,\sigma_i\sigma_{i+1},\dots,\sigma_{n+1}) \\
   && +(-1)^{n+1}L(\sigma_{n+1},X_0)c(\sigma_1,\dots,\sigma_n).
\end{eqnarray*}
For $n=0$ and $\sigma\colon X\r Y$, $d(c)(\sigma)=L(1,\sigma)c(X)-L(\sigma,1)c(Y)$. 

The cohomology of $\C{C}$ with coefficients in $L$ is the
cohomology of the complex $F^*(\C{C},L)$. It is denoted by
$H^*(\C{C},L)$. For the functorial behaviour of cohomology of categories see \cite{fcc}.

If the bimodule $L$ takes values in the category of $k$-modules for $k$ a commutative ring then $H^*(\C{C},L)$ is a graded $k$-module.
\end{defn}

It follows from the very definition that given an additive category $\C{A}$ the cohomology group $H^0(\C{A},\hom_\C{A})$ is isomorphic to the endomorphism group of the identity functor in $\C{A}$, see \cite[Proposition 3.2]{cat}. Given an $R$-module $M$ the following cohomology groups of $\C{P}(R)$ coincide with Mac Lane cohomology by \cite[Theorem A and Corollary 3.11]{cat}.
\begin{eqnarray}\label{jp}
H^*(\C{P}(R),\hom_R(-,M\otimes_R -))&\cong& HML^*(R,M).
\end{eqnarray}
For $R=M=\Z/2$ the $0$-dimensional group was computed in \cite{mlcff}, $H^0(\Z/2,\Z/2)\cong \Z/2$. From this isomorphism we derive the following proposition.

\begin{prop}\label{solo1} 
The identity functor in $\C{P}(\Z/2)$ has only two endomorphisms, the identity and the zero endomorphisms. In particular it has only one automorphism, the identity morphism.
\end{prop}

Let us recall now the Pirashvili-Waldhausen homology of categories defined in \cite{mlhthh}, which is defined in a dual way to Baues-Wirsching cohomology.

\begin{defn}\label{hc}
Let $\C{C}$ be a category and let $L$ be a $\C{C}$-bimodule.
The \emph{Pirashvili-Waldhausen
complex} $F_*(\C{C},L)$ is given by the direct sums
\begin{eqnarray*}
F_n(\C{C},L)&=&\bigoplus_{X_0\st{\sigma_1}\l\cdots\st{\sigma_n}\l X_n}L(X_0,X_n),\;\;n\geq0.
\end{eqnarray*}
An element $c\in L(X_0,X_n)$ regarded as an $n$-chain in the direct factor correspoding to $X_0\st{\sigma_1}\l\cdots\st{\sigma_n}\l X_n$ will be denoted by
$c\cdot(\sigma_1,\dots,\sigma_n)$. The differential is defined by
\begin{eqnarray*}
  d(c\cdot(\sigma_1,\dots,\sigma_n)) &=& L(\sigma_1,X_n)(c)\cdot(\sigma_2,\dots,\sigma_n)\\
   && +\sum_{i=1}^{n-1}(-1)^ic\cdot(\sigma_1,\dots,\sigma_i\sigma_{i+1},\dots,\sigma_n) \\
   && +(-1)^nL(X_0,\sigma_n)(c)\cdot(\sigma_1,\dots,\sigma_{n-1}),
\end{eqnarray*}
for $n\geq 2$, and for $n=1$ and $\sigma\colon X\r Y$ we have
$d(c\cdot\sigma)=L(\sigma,X)(c)\cdot X-L(Y,\sigma)(c)\cdot Y$. This differential differs from the one defined in \cite{mlhthh} by a sign. Of course this does not affect the homology.

The homology of $\C{C}$ with coefficients in $L$ is the
homology of the complex $F_*(\C{C},L)$. It is denoted by
$H_*(\C{C},L)$. 

If $L$ takes values in the category of $k$-modules for $k$ a commutative ring then $H_*(\C{C},L)$ is a graded $k$-module.
\end{defn}

By \cite{mlhthh} the homology of $\C{P}(R)$ is naturally isomorphic to the Mac Lane homology of the ring $R$ and to the topological Hochschild homology of the Eilenberg-Mac Lane ring spectrum $HR$ in the following case,
\begin{eqnarray}
\label{iso2} H_*(\C{P}(R),\hom_R)&\cong& HML_*(R)\\
\nonumber &\cong& THH_*(HR).
\end{eqnarray}

In the following theorem we construct a \emph{universal coefficients spectral sequence} for the Mac Lane cohomology of commutative rings which is crucial for the computations in this paper.

\begin{thm}
For any commutative ring $R$ and any $R$-module $M$ there is a spectral sequence
\begin{eqnarray*}
E^{p,q}_2=\ext^p_R(HML_q(R),M)&\Longrightarrow&HML^{p+q}(R,M).
\end{eqnarray*}
Here $M$ is regarded as a symmetric $R$-bimodule for the definition of the Mac Lane cohomology.
\end{thm}

\begin{proof}
Given two finitely generated projective $R$-modules $P, Q$ there are natural isomorphisms
$$\hom_R(\hom_R(Q,P),M)\st{\alpha}\longleftarrow M\otimes_R\hom_R(P,Q)\st{\beta}\To\hom_R(P,M\otimes_R Q),$$
defined by
\begin{eqnarray*}
\alpha (m\otimes f)(g)&=&m\cdot\trace(fg),\;\; g\in \hom_R(Q,P),\\
\beta(m\otimes f)(x)&=&m\otimes f(x),\;\;x\in P.
\end{eqnarray*}
Indeed the homomorphisms $\alpha$ and $\beta$ are clearly natural and biadditive, so it is enough to check that they are isomorphisms for $P=Q=R$, and this case is trivial.

The isomorphisms $\alpha$ and $\beta$ determine a cochain isomorphism
\begin{eqnarray*}
\hom_R(F_*(\C{P}(R),\hom_R),M)&\cong&F^*(\C{P}(R),\hom_R(-,M\otimes_R-)),
\end{eqnarray*}
compare Definitions \ref{cc} and \ref{hc}.

Now if $E^*$ is an injective resolution of $M$ then the spectral sequence of the statement is the spectral sequence of the bicomplex
$$\hom_R(F_*(\C{P}(R),\hom_R),E^*).$$
\end{proof}

The following proposition is a simple application of the previous theorem.

\begin{prop}\label{34}
For any $\Z/4$-module $M$ there is a natural isomorphism
\begin{eqnarray*}
HML^3(\Z/4,M)&\cong&\hom_{\Z/4}(\Z/2,M).
\end{eqnarray*}
\end{prop}

\begin{proof}
The topological Hochschild homology of $\Z/4$ is computed in \cite{thhzpn}. The lower homology groups are
\begin{equation*}
HML_n(\Z/4)\;\;\cong\;\;\left\{\begin{array}{ll}
\Z/4,&n=0,\\
0,&n=1,\\
\Z/4,&n=2,\\
\Z/2,&n=3.
\end{array}\right.
\end{equation*}
This implies that the $E_2$-term of the universal coefficients spectral sequence satisfies $E^{p,q}_2=0$ in case $p>0$ and $q<3$. Therefore the isomorphism of the statement follows. 
\end{proof}

\section{Toda brackets in triangulated and homotopy categories}

We recall from \cite[13]{shc} the definition of Toda bracket in a triangulated category.

\begin{defn}\label{tbt}
Let $\C{T}$ be a triangulated category with translation functor $\Sigma$. Given three composable morphisms in $\C{T}$ 
$$W\st{f}\To X\st{g}\To Y\st{h}\To Z,$$ 
such that $gf=0$ and $hg=0$, the \emph{Toda bracket}
$$\grupo{h,g,f}\in\frac{\C{T}(\Sigma W,Z)}{h\cdot\C{T}(\Sigma W,Y)+\C{T}(\Sigma X,Z)\cdot(\Sigma f)}$$
is defined in the following way. If 
$$W\st{f}\To X\st{i}\To C_f\st{q}\To \Sigma W$$
is an exact triangle in $\C{T}$ the axioms of a triangulated category imply the existence of a commutative diagram
$$\xymatrix{W\ar[r]^f\ar@{=}[d]&X\ar[r]^i\ar@{=}[d]&C_f\ar[r]^q\ar[d]^a&\Sigma W\ar[d]^b\\
W\ar[r]_f&X\ar[r]_g&Y\ar[r]_h&Z}$$
The morphism $b$ is a representative of the Toda bracket $\grupo{h,g,f}$. 
The Toda bracket does not depend on the choices made for its definition since we divide out the indeterminacy. Moreover,
all representatives of $\grupo{h,g,f}$ can be obtained in this way.
\end{defn}

\begin{exm}\label{exefu}
We derive from the very definition that the Toda bracket of the morphisms of an exact triangle like (\ref{tri}) contains the identity morphism $1_{\S A}\in\grupo{q,i,f}$. In fact exact triangles are characterized among candidate triangles by this property and a further homological condition, see \cite[Therorem 13.2]{shc}. This and Proposition \ref{lose} imply that, in the triangulated category $\C{P}(\Z/4)$ of Theorem \ref{rara}, the sequence of morphisms $$\Z/4\st{2}\To\Z/4\st{2}\To\Z/4\st{2}\To\Z/4$$ has Toda bracket
$$0 \;\neq \;\grupo{2,2,2} \;\in \;\frac{\C{P}(\Z/4)(\Z/4,\Z/4)}{2\cdot\C{P}(\Z/4)(\Z/4,\Z/4)+\C{P}(\Z/4)(\Z/4,\Z/4)\cdot 2} \;\cong\; \Z/2.$$
\end{exm}

Toda brackets in the homotopy category of a stable model category are determined by a class in Baues-Wirsching cohomology of categories. 

A $3$-dimensional class in the Baues-Wirsching cohomology of a category with zero object defines Toda brackets as we now indicate. The \emph{Toda category} $\mathbf{Toda}$ introduced in \cite{ccglg} has five objects $*,1,2,3,4$. The object $*$ is a zero object, and there are only three morphisms which are neither identities nor zero morphisms, namely $j_n\colon n\r n+1$ for $n=1,2,3$. The Toda category can be represented by the following commutative diagram
$$\mathbf{Toda}\;\;=\;\;\xymatrix{1 \ar[r]|{j_1}\ar@/^20pt/[rr]^0&2\ar[r]|{j_2}\ar@/_20pt/[rr]_0&3\ar[r]|{j_3}&4.}$$
If $M$ is a $\mathbf{Toda}$-bimodule with $M(*,-)=0$ and $M(-,*)=0$ then
\begin{eqnarray*}
H^3(\mathbf{Toda},M)&\cong&\frac{M(1,4)}{M(1,j_3)M(1,3)+M(j_1,4)M(2,4)}.
\end{eqnarray*}
The isomorphism follows from the fact that, with this kind of coefficients, the Baues-Wirsching cohomology can be computed by using cochains which are normalized with respect to identity and zero morphisms, see \cite{ccglg}. This isomorphism sends a cohomology class represented by such a normalized $3$-cocycle $c$ to the class of $c(j_3,j_2,j_1)$. 

\begin{defn}
Given a category $\C{C}$ with zero object $*$, a $\C{C}$-bimodule $L$ with $L(*,-)=0$ and $L(-,*)=0$, and a cohomology class $\zeta\in H^3(\C{C},L)$, then the \emph{$\zeta$-Toda bracket} of a diagram in $\C{C}$
$$W\st{f}\To X\st{g}\To Y\st{h}\To Z,$$ 
with $gf=0$ and $hg=0$, is defined as follows. The diagram corresponds to a functor $\varphi\colon\mathbf{Toda}\r\C{C}$ with $\varphi(*)=*$, $\varphi(j_1)=f$, $\varphi(j_2)=g$, and $\varphi(j_2)=h$, and the $\zeta$-Toda bracket is the pull-back along $\varphi$ of the cohomology class $\zeta$
\begin{equation*}
\grupo{h,g,f}_\zeta\;=\;\varphi^*\zeta \;\in\; H^3(\mathbf{Toda},L(\varphi,\varphi))\;\cong\;\frac{L(W,Z)}{L(W,h)L(W,Y)+
L(f,Z)L(X,Z)}.
\end{equation*}
If $c$ is a cocycle representing $\zeta$ which is normalized with respect to identity and zero morphisms then $\grupo{h,g,f}_\zeta$ is represented by $c(h,g,f)$.
\end{defn}

If a triangualted category $\C{T}=\ho\C{M}$ is the homotopy category of a stable model category $\C{M}$ then Toda brackets in $\C{T}$ are determined by the \emph{universal Toda bracket}
$$\grupo{\C{M}}\in H^3(\C{T},\hom_{\C{T}}(\Sigma,-)),$$
i.e. Toda brackets in the triangulated category $\ho\C{M}$ are $\grupo{\C{M}}$-Toda brackets in the sense of the previous definition. This universal Toda bracket is defined in \cite{ccglg} on the homotopy category of fibrant-cofibrant objects, but it can be placed in the cohomology of the full homotopy category by using the fact that equivalences of categories induce isomorphisms in Baues-Wirsching cohomology. By \cite[Theorem 13.2]{shc} the universal Toda bracket also determines the exact triangles in $\ho\C{M}$.

\begin{prop}\label{0es}
For any $\zeta\in H^3(\C{P}(\Z/4),\hom_{\Z/4})\cong\Z/2$ the $\zeta$-Toda bracket of $$\Z/4\st{2}\To\Z/4\st{2}\To\Z/4\st{2}\To\Z/4$$ vanishes.
\end{prop}

\begin{proof}
The isomorphism $H^3(\C{P}(\Z/4),\hom_{\Z/4})\cong\Z/2$ follows from the isomorphism (\ref{jp}) and Proposition \ref{34}. Actually the isomorphism class of this group is not important for the proof. We just want to remark that if it were trivial there would be nothing to prove, but this is not the case, so the argument below is needed.

Let $\varphi\colon\mathbf{Toda}\r\C{P}(\Z/4)$ be the functor corresponding to $$\Z/4\st{2}\To\Z/4\st{2}\To\Z/4\st{2}\To\Z/4,$$ and let $i\colon\Z/2\hookrightarrow\Z/4$ be the inclusion. The functorial behaviour of Baues-Wirsching cohomology implies that the following diagram commutes
\begin{equation*}\tag{a}
\xymatrix@C=60pt{H^3(\C{P}(\Z/4),\hom_{\Z/4}(-,\Z/2\otimes_{\Z/4}-))\ar[r]^-{\hom_{\Z/4}(-,i\otimes-)_*}\ar[d]_{\varphi^*}&H^3(\C{P}(\Z/4),\hom_{\Z/4})
\ar[d]^{\varphi^*}\\
H^3(\mathbf{Toda},\hom_{\Z/4}(\varphi,\Z/2\otimes_{\Z/4}\varphi))\ar[r]_-{\hom_{\Z/4}(\varphi,i\otimes\varphi)_*}&H^3(\mathbf{Toda},\hom_{\Z/4}(\varphi,\varphi))}
\end{equation*}
Here an upper $*$ indicates the pull-back along a functor and the lower $*$ indicates a push-forward along a bimodule morphism.

The isomorphism in (\ref{jp}) together with Proposition \ref{34} prove that we have a commutative diagram
\begin{equation*}
\xymatrix@C=60pt{H^3(\C{P}(\Z/4),\hom_{\Z/4}(-,\Z/2\otimes_{\Z/4}-))\ar[r]^-{\hom_{\Z/4}(-,i\otimes-)_*}\ar[d]_\cong&H^3(\C{P}(\Z/4),\hom_{\Z/4})
\ar[d]^\cong\\
\hom_{\Z/4}(\Z/2,\Z/2)\ar[r]^-{\hom_{\Z/4}(\Z/2,i)}\ar[d]_\cong&\hom_{\Z/4}(\Z/2,\Z/4)\ar[d]^\cong\\
\Z/2\ar@{=}[r]&\Z/2}
\end{equation*}
Hence the upper horizontal homomorphism in (a) is an isomorphism.

On the other hand we have another commutative diagram
\begin{equation*}
\xymatrix@C=60pt{H^3(\mathbf{Toda},\hom_{\Z/4}(\varphi,\Z/2\otimes_{\Z/4}\varphi))\ar[r]^-{\hom_{\Z/4}(\varphi,i\otimes\varphi)_*}\ar[d]_\cong&
H^3(\mathbf{Toda},\hom_{\Z/4}(\varphi,\varphi))\ar[d]^\cong\\
\frac{\hom_{\Z/4}(\Z/4,\Z/2)}{2\cdot \hom_{\Z/4}(\Z/4,\Z/2)+\hom_{\Z/4}(\Z/4,\Z/2)\cdot 2}\ar[r]^-{\hom_{\Z/4}(\Z/4,i)}\ar[d]_\cong&
\frac{\hom_{\Z/4}(\Z/4,\Z/4)}{2\cdot \hom_{\Z/4}(\Z/4,\Z/4)+\hom_{\Z/4}(\Z/4,\Z/4)\cdot 2}\ar[d]^\cong\\
\hom_{\Z/4}(\Z/4,\Z/2)\ar[r]^-{\hom_{\Z/4}(\Z/4,i)}\ar[d]_\cong&\hom_{\Z/4}(\Z/4,\Z/4)/(2\cdot 1_{\Z/4})\ar[d]^\cong\\
\Z/2\ar[r]_0&\Z/2}
\end{equation*}
Therefore the lower horizontal homomorphism in diagram (a) is $0$, so the vertical homomorphism in the right of diagram (a) must also be zero. This finishes the proof.
\end{proof}

Now we are ready to prove Theorem \ref{raro}.

\begin{proof}[Proof of Theorem \ref{raro}]
Suppose by the contrary that there is a fully faithful exact functor $\psi\colon\C{P}(\Z/4)\r\ho\C{M}$. Such a functor preserves Toda brackets, therefore Toda brackets in $\C{P}(\Z/4)$ are $\psi^*\grupo{\C{M}}$-Toda brackets, hence 
by Proposition \ref{0es} the Toda bracket of $$\Z/4\st{2}\To\Z/4\st{2}\To\Z/4\st{2}\To\Z/4$$ vanishes. This contradicts the computation in Example \ref{exefu}.
\end{proof}

\bibliographystyle{amsalpha}
\bibliography{Fernando}
\end{document}